\documentclass[12pt]{article}
\usepackage[all]{xy}
\usepackage{amsmath}
\usepackage{amssymb}
\usepackage{amsthm}

\newtheorem{theorem}{Theorem}
\newtheorem{lemma}[theorem]{Lemma}
\newtheorem{proposition}[theorem]{Proposition}
\newtheorem{corollary}[theorem]{Corollary}

\newtheorem{fact}[theorem]{Fact}
\theoremstyle{definition}
\newtheorem{definition}[theorem]{Definition}

\begin{document}
\newcommand{\xx}{\mathbf{x}}
\newcommand{\yy}{\mathbf{y}}
\newcommand{\conj}{\mathbf{conj}}
\newcommand{\hypfails}{\mathbf{hypfails}}
\newcommand{\concholds}{\mathbf{concholds}}
\newcommand{\inl}{\mathbf{inl}}
\newcommand{\inr}{\mathbf{inr}}
\newcommand{\inlinr}{\mathbf{inlinr}}
\newcommand{\reasons}{\mathrm{Reasons}}
\newcommand{\reasonsinf}{\mathrm{Reasons}^\infty}
\newcommand{\rank}{\mathrm{rank}}
\newcommand{\N}{\mathbb{N}}
\newcommand{\0}{\mathbf{0}}
\renewcommand{\subset}{\subseteq}
\newcommand{\straightforward}{\begin{proof} Straightforward.\end{proof}}
\newcommand{\lev}{\mathrm{Lev}}
\newcommand{\num}{\#}
\newcommand{\ran}{\mathrm{ran}}
\renewcommand{\implies}{\rightarrow}
\newcommand{\meet}{\wedge}	
\newcommand{\powerset}{\mathop{\mathcal{P}}}
\bibliographystyle{amsplain}
\title{Embeddings into Free Heyting Algebras and Translations into Intuitionistic Propositional Logic}
\author{Michael O'Connor\\
          Department of Mathematics\\ Cornell University\\ Ithaca NY 14853\\
          \texttt{oconnor@math.cornell.edu}}
\maketitle
\begin{abstract}
We find a translation with particularly nice properties from intuitionistic propositional logic
in countably many variables to intuitionistic propositional logic in two variables. In addition, the existence
of a possibly-not-as-nice translation from any countable logic into intuitionistic propositional logic in two variables is shown.  The nonexistence of a translation from classical logic into intuitionistic propositional logic
which preserves $\wedge$ and $\vee$ but not necessarily $\top$ is proven.  These results about translations follow from additional results about embeddings into free Heyting algebras.
\end{abstract}

\section{Introduction}
Intuitionistic logic has been explored for many years as a
language for computer science, with a guiding principle being the Brouwer-Heyting-Kolmogorov
interpretation, under which intuitionistic proofs of implication 
are functions and existence proofs require witnesses. Higher-order intuitionistic systems
which can express a great deal of mathematics, such as Girard's System
F and Martin-L\"of's type theory (good references are~\cite{girardpat} and~\cite{beeson}), have been developed and implemented by
prominent computer scientists such as Constable, Huet and Coquand (see~\cite{nuprl} and~\cite{coq}). With
all this development and with the existence of well-established topological, Kripke, and
categorical semantics for intuitionistic systems, it may come as a
surprise that many fundamental structural properties of intuitionistic
propositional calculus have not been developed.
By way of contrast, corresponding issues for classical logics have been settled for
at least 75 years. 

Heyting algebras are an equationally defined class of algebras with operations 
$\vee$, $\meet$, and $\rightarrow$ and constants $\bot$ and $\top$ (representing ``or,'' ``and,'' ``implies,'' ``false,'' and ``true'' respectively) that stand in the same relation to intuitionistic propositional
logic that Boolean algebras do to classical propositional logic. What follows is a very brief introduction to free Heyting
algebras and a summary of the results that will be presented in this paper.

\newcommand{\iv}{\vdash_{\mathcal{I}}}
\newcommand{\cv}{\vdash_{\mathcal{C}}}

For each $n\in\N$, let $V_n = \{x_1,\ldots, x_n\}$ and let $F_n$ be the set of
propositional sentences in variables $V_n$. 
Let $\simeq^n_i$ and $\simeq^n_c$ be the intuitionistic and classical logical equivalence
relations respectively.  

The classical Lindenbaum algebra $B_n$ is then defined as $F_n/\simeq^n_c$ 
and the intuitionistic Lindenbaum algebra $H_n$ is defined as $F_n/\simeq^n_i$.
The operations $\meet$ and $\vee$ and the constants $\top$ and $\bot$ are naturally defined on $B_n$ and the operations
$\meet$, $\vee$, and $\rightarrow$ and the constants $\top$ and $\bot$ are naturally defined on $H_n$. $B_n$
is then isomorphic to the free Boolean algebra on $n$ generators and $H_n$ is the free Heyting algebra
on $n$ generators. As usual, the order $\leq$ may be defined from $\wedge$ (or from $\vee$). Like all Heyting algebras, each $H_n$ is also a distributive lattice.

The analogous statements are true for $V_\omega = \{x_1,x_2,\ldots\}$, $F_\omega$, $B_\omega$, and
$H_\omega$.
%
%
%

The structure of each $B_n$ and of $B_\omega$ is well understood. However, among the free Heyting algebras,
only $H_1$ is completely understood. It is known from~\cite{nishimura} that if we let $\phi_1 = \neg x_1$,
$\psi_1 = x_1$, $\phi_{i+1} = \phi_i \rightarrow \psi_i$, and $\psi_{i+1} = \phi_i \vee \psi_i$, then
each propositional formula in the single variable $x_1$ is intuitionistically equivalent to exactly one formula
in $\{\bot\}\cup  \{\phi_i\mid i \in\omega\} \cup  \{\psi_i\mid i \in\omega\}\cup\{\top\}$.
Further, we can easily write down conditions characterizing the order on those formulas, so that the structure of $H_1$ is completely characterized.

Although the structure of $H_n$ for $n\geq 2$ and of $H_\omega$ is not fully understood, there are a number of facts known.  Although not a complete list, the reader is referred to~\cite{bellissima}, \cite{ghilardi}, \cite{butz},
and \cite{darniere}.  A very useful construction is contained in~\cite{bellissima} which will we avail ourselves of in this paper and which is described in Section~\ref{notation} below.

%
%


The results of this paper are as follows: 
%
There is a lattice-embedding from $H_\omega$ into $H_2$.
This obviously implies that there is a lattice-embedding from $H_m$ into $H_n$
for any $m \geq 1$, $n\geq 2$, but in these cases, more is true:
For any $n \geq 2$ and $m\geq 1$,
there is a $\phi$ and a 
$\psi\in H_n$ such that $[\phi,\psi] := \{\rho \in H_n \mid \phi \leq \rho \leq \psi\}$
and $H_m$ are isomorphic as lattices.  In addition, the isomorphism from $[\phi,\psi]$
to $H_m$ can be extended to all of $H_n$, so that there is a surjective lattice-homomorphism
from $H_n$ to $H_m$.

Furthermore, we will show that any countable partial order can be order-embedded into $H_2$,
and that the countable atomless boolean algebra $B_\omega$ cannot be lattice-embedded into $H_\omega$.

Some of these results also have significance in terms of \emph{translations} into intuitionistic propositional
logic, a notion which we now define.
%


Let $\vdash$ be the intuitionistic consequence relation.  Define a consequence-respecting translation from $n$-variable intuitionistic logic into $m$-variable intuitionistic logic
to be a function $f\colon F_n\to F_m$ such that for all $\emptyset\ne\Gamma\subset F_n$,
$\phi\in F_n$, $\Gamma \vdash \phi$ iff $f(\Gamma) \vdash f(\phi)$.


Define a tautology-respecting translation from $n$-variable intuitionistic logic into $m$-variable
intuitionistic logic to be a function $f\colon F_n \to F_m$ such that for all $\phi\in F_n$,
$\vdash \phi$ iff $\vdash f(\phi)$.

The term ``respecting'' is used to emphasize that the property of being consequence-respecting and
the property of being tautology-respecting are stronger than the property of being consequence-preserving
and tautology-preserving respectively.


We define a $(\wedge, \vee)$-preserving translation from $n$-variable propositional logic
to $m$-variable propositional logic to be a function $f\colon F_n\to F_m$ such that 
for all $\phi$, $\psi\in F_n$, $f(\phi\wedge\psi) = f(\phi)\wedge f(\psi)$ and
$f(\phi\vee\psi) = f(\phi)\vee f(\psi)$.


We make the obvious modifications to the definitions for translations from classical logic to intuitionistic logic
and for $\omega$-variable logics.

Thus, G\"odel's double-negation translation (see~\cite{nerodelambda} or~\cite{beeson}) is a tautology-respecting but not consequence-respecting
or $(\wedge,\vee)$-preserving translation from $\omega$-variable classical logic to $\omega$-variable
intuitionistic logic and Gentzen's translation (again, see~\cite{nerodelambda} or~\cite{beeson}) is a tautology-respecting and consequence-respecting
but not $(\wedge,\vee)$-preserving translation from $\omega$-variable classical logic to $\omega$-variable
intuitionistic logic.  Both of these translations may be restricted to be from $n$-variable classical
logic to $n$-variable intuitionistic logic for any $n$.


Some of the results of this paper may then be restated as follows: There is a consequence- and tautology-respecting, 
$(\wedge,\vee)$-preserving translation of $\omega$-variable intuitionistic logic into 2-variable 
intuitionistic logic.
We also get consequence- and tautology-respecting translations (which aren't $(\wedge,\vee)$-preserving) of $\omega$-variable classical logic into 2-variable intuitionistic logic by composing with
Gentzen's translation.  This translation may be read off explicitly from the proof contained in this paper together
with the construction of~\cite{bellissima}.

The disjunction property of intuitionistic logic implies that there can be no tautology-respecting translation of classical logic into intuitionistic logic.  In addition, we will show that there is no merely consequence-respecting, $(\wedge,\vee)$-preserving translation of $\omega$-variable
classical logic into $\omega$-variable intuitionistic logic (and thus not into
$n$-variable intuitionistic logic for any $n$).

The result that any countable partial order can be embedded in $H_2$ implies that \emph{any} logic
may be translated in a consequence- and tautology-respecting but not necessarily $(\wedge,\vee)$-preserving way into 2-variable intuitionistic logic, as long as the logic is countable.

The author would like to acknowledge Richard Shore and Anil Nerode for many useful conversations and specific comments on this paper, and his parents for their love and guidance.
\section{Notation, Terminology, and Bellissima's Construction}\label{notation}
As above, let  $V_n = \{x_1,x_2,\ldots, x_n\}$. 

We will use Bellissima's construction (\cite{bellissima}) of, for each $n$, a Kripke model $K_n$ over $V_n$ satisfying Propositions~\ref{exact} and~\ref{phidef} below.  The construction and relevant facts about it will be stated here.



\newcommand{\nodes}{\mathop{\text{Nodes}}}
\newcommand{\kleq}{\mathop{\leq}}
\newcommand{\node}{\mathop{\text{node}}}

Given a Kripke model $K$ over $V_n$, let $\nodes(K)$ be the set of nodes of $K$ and let $\kleq(K)$ be the (non-strict) partial order
on $\nodes(K)$ given by $K$. 
If no confusion will result, we may use $K$ in place of $\nodes(K)$. Given $\alpha\in K$, let
$w(\alpha) = \{x_i\in V_n\mid \alpha\Vdash x_i\}$.

We will define a Kripke model $K_n$ in stages, so that $K_n = \bigcup_i K^i_n$ where the $K^i_n$ are defined as follows:

$\nodes(K^0_n) = \powerset(V_n)$ and $\kleq(K^0_n) = \{(\alpha,\alpha) \mid \alpha\in \nodes(K^0_n)\}$.
For clarity, when we want to emphasize that we are thinking of $U\subset V_n$ as
a node, we may write $\node(U)$ or $\node_{K_n}(U)$. If we want to also note
that it is in $K^0_n$, we may write $\node^0(U)$ or $\node^0_{K_n}(U)$. For $U\subset V_n$, we let $w(\node(U)) = U$.

Given $K^0_n,\ldots, K^i_n$, let $\mathcal{T}_{i+1}$ be the set of subsets $T$ of $\bigcup_{j = 0}^i K^j_n$
such that $T\cap (K^i_n- K^{i-1}_n) \ne \emptyset$ and such that the elements of $T$ are pairwise incomparable with respect
to $\kleq(K^i_n)$.  Then we define $\nodes(K^{i+1}_n) - \nodes(K^i_n)$ to be
\[ 
\{ \langle T,U\rangle \mid T\in\mathcal{T}_{i+1}, U\subset V_n, U\subset\bigcap_{\alpha\in T} w(\alpha)
\text{ and if }T = \{\beta\},\text{ then }U\subsetneq w(\beta)\}
\]
For clarity, when we want to emphasize that we are thinking of $\langle T,U\rangle$ as a node,
we may write $\node(\langle T,U\rangle)$ or $\node_{K_n}(\langle T,U\rangle)$.
If we want to also note that it is in $K^{i+1}_n$,
we may write $\node^{i+1}(\langle T,U\rangle)$ or $\node^{i+1}_{K_n}(\langle T,U\rangle)$.

We declare that $w(\langle T,U\rangle) = U$ and we
let $\kleq(K^{i+1}_n)$ be the reflexive transitive closure of
$\kleq(K^i_n)\cup \{(\langle T,U\rangle,\beta) \mid \langle T,U\rangle \in K^{i+1}_n - K^i_n, \beta\in T\}$.


Let $k(\phi)$ denote the set of nodes in $K_n$ which force $\phi$, for $\phi$ a propositional formula in
$n$ variables.

\begin{proposition}[\cite{bellissima}]\label{exact} For $\phi$ and $\psi$ propositional formulas in $x_1,\ldots, x_n$,
$\phi\vdash \psi$ iff $k(\phi)\subset k(\psi)$.
\end{proposition}

\begin{proposition}[\cite{bellissima}]\label{phidef} For each node $\alpha\in K_n$, there is a $\phi_\alpha$ such that
$k(\phi_\alpha) = \{\beta\in K_n\mid \beta \geq \alpha\}$ and there is a $\phi'_\alpha$ such that 
$k(\phi'_\alpha) = \{\beta\in K_n\mid \beta\not\leq \alpha\}$.
\end{proposition}

%

%
%

We now fix some terminology.

If $\alpha < \beta$ are nodes in some Kripke model, then $\beta$ is called a successor
of $\alpha$ and $\alpha$ a predecessor of $\beta$.  If there is no $\gamma$
with $\alpha < \gamma< \beta$, then $\beta$ is called an immediate successor of 
$\alpha$. The assertions ``$\alpha$ is above $\beta$'' and ``$\beta$ is below $\alpha$''
both mean $\alpha \geq \beta$.

For any node $\alpha$ in any Kripke model, $s(\alpha)$ is the set of $\alpha$'s immediate
successors.

If $\alpha\in K_n$, then $\phi_\alpha$ and $\phi'_\alpha$ are as in Proposition~\ref{phidef}.

For each $m$, $\lev^n_m = K^m_n - K^{m-1}_n$. This may also be called
$\lev_m$ if $n$ is clear from context and may be denoted in English as ``level $m$.''

%

If $\alpha\in K_n$, then $\lev(\alpha)$ is the unique $i$ such that $\alpha\in \lev^n_i$.
Note that if $\alpha \leq \beta$, $\lev(\alpha) \geq \lev(\beta)$.
%

If $T$ is a set of nodes in $K_n$, let $r(T) = \{\alpha\in T\mid \neg(\exists \beta\in T)\,
(\beta<\alpha)\}$.  Thus, for example, for any $T$ with $|T|\geq 2$, 
$\langle r(T),\emptyset\rangle \in K_n$. 
The following facts will be used below and follow without much difficulty directly from the construction.
\begin{fact}\label{biglevel} For $n\geq 2$ and $m\geq 0$, $|\lev^n_{m+1}| > |\lev^n_{m}|$. In particular, there are arbitrarily large levels of
$K_n$.
\end{fact}
The following fact is a more general version of the preceding fact.
\begin{fact}\label{bignotnot} Let $S\subset K_n$, $|S| \geq 3$ and let each element of $S$ be at the same level. Let $S'$
be the downward closure of $S$. Then $|S'\cap\lev^n_{m+1}| > |S'\cap\lev^n_m|$ for any $m$ greater than or equal
to the common level of the elements of $S$	. 
\end{fact}

\section{A Lattice Embedding from $H_m$ to $H_n$ for $m \geq 1$, $n \geq 2$}

\begin{theorem} Let $n\geq 2$, $m\geq 1$. Then
there are $\phi$, $\psi \in H_n$ such that $H_m$ is isomorphic to $[\phi,\psi]$.
In addition, the isomorphism from $[\phi,\psi]$ to $H_m$ can be extended to a surjective lattice-homomorphism
from $H_n$ to $H_m$.
\end{theorem}
\begin{proof}
The main work is contained in the following proposition.

\begin{proposition}\label{embnoint} Let $m\geq 2$ and $n$ be such that there is a level $\lev^n_i$ of $K_n$ 
and a set $A\subset\lev^n_i$ such that $|A| = m$ and each $\alpha\in A$ has
some immediate successor not above any other $\alpha'\in A$.  Then there is a $\phi$, $\psi\in H_n$ such
that $H_m$ is lattice-isomorphic to $[\phi,\psi]$ and in addition, the isomorphism from $[\phi,\psi]$ to $H_m$ can be extended to a surjective lattice-homomorphism
from $H_n$ to $H_m$.
\end{proposition}
\begin{proof}
Fix $A$ and $i$ from the hypothesis. 

Let $A=\{\alpha_1,\ldots,\alpha_m\}$.  Let $\phi$ be
\[ \bigvee_i \phi_{\alpha_i}. \]

For each $A'\subset A$, let $\gamma_{A'}$ be 
the node $\langle r(T),\emptyset\rangle$ where $T = A'\cup \bigcup_{\alpha\notin A'} s(\alpha)$. This is valid
as the elements of $r(T)$ are pairwise incomparable and $|r(T)| \geq 2$ since $m \geq 2$.

Note that $\gamma_{A'}$ is at level $i+1$ if $A'$ is nonempty and at level $i$ if $A'$ is empty.
Since each $\alpha_i$ has a successor not above any other $\alpha_j$, 
if $A' \ne A''$, $\gamma_{A'}\ne \gamma_{A''}$.

Let $S = \{\rho\in K^{i+1}_n\mid (\forall A'\subset A)\,(\rho\not\geq\gamma_{A'})\}$
and let $\psi$ be $\psi_0\wedge \psi_1$ where $\psi_0$ is
\[ \left[\neg\neg(\bigvee_{A'\subset A}\phi_{\gamma_{A'}})\right]\]
and $\psi_1$ is
\[ \bigwedge_{\rho\in S} \phi'_\rho\]

Define a function $g$ with domain $K_m$ as follows:

1. $g(\node_{K_m}(U)) = \gamma_{A'}$ where $A' = \{\alpha_k \mid x_k\in V_m - U\}$.

2. $g(\node_{K_m}^{j+1}(\langle T,U\rangle)) = \langle r(T'),\emptyset\rangle$, where 
$T' = \{g(\delta) \mid \delta\in T\}\cup \{\alpha_k \mid x_k \in V_m - U\}$.

We will show that the range of $g$ is contained in $K_n$.  By induction, what we must show is
that $\langle r(T'),\emptyset\rangle$ is in $K_n$, which will hold as long as $|r(T')| \geq 2$.

\begin{lemma}\label{injection} The function $g$ is into $K_n$ and preserves order and nonorder.
For all $\beta\in K_m$ and $x_k\in V_m$, $\beta \Vdash x_k$ iff $g(\beta)\not\leq \alpha_k$.
\end{lemma}
\begin{proof}
We will prove by induction on $i$ that $g$ restricted to $K^i_m$ satisfies the conditions in the statement of the
lemma.

%


For $i = 0$, observe that $\{g(\node^0(U))\mid U\subset V_m\}$ is pairwise incomparable and
that if $U\ne U'$, $g(\node^0(U))\ne g(\node^0(U'))$ as they have different immediate successors. It is also
the case that for all $\node^0(U)\in K_m^0$ and $x_k\in V_m$, $\node^0(U) \Vdash x_k$ iff
$x_k \in U$ iff $\gamma_{A'}\not\leq \alpha_k$, where $A' = \{\alpha_k \mid x_k \in V_m - U\}$.

Finally, since each $\gamma_{A'}$ is in $K_n$, the range of $g$ restricted to $K^0_m$ is contained in $K_n$.

Now suppose $g$ restricted to $K^i_m$ satisfies the hypotheses in the statement of the lemma.


We first show that the range of $g$ restricted to $K^{i+1}_m$
is contained in $K_n$. Let $\langle T,U\rangle \in \lev^m_{i+1}$. If $|T| \geq 2$, then $|r(T)| \geq 2$
and we are done. If $|T| = \{\beta\}$, then $U \subsetneq w(\beta)$ and $T'$ must contain both
$g(\beta)$ and $\alpha_k$, where $x_k \in w(\beta) - U$.  Since $\beta \Vdash x_k$, $g(\beta)\not \leq
\alpha_k$.  Since it is fairly easy to see that each $\alpha_k$ is not less than any element of the range of
$g$, we must have $|r(T')| \geq 2$.

It is immediate then that $g$ restricted to $K^{i+1}_m$ is preserves order and the immediate successor
relation. Each element
of $\lev^m_{i+1}$ is of the form $\node^{i+1}(T,U)$. Observe that if $U\ne U'$
and $\langle T, U\rangle,
\langle T,U'\rangle\in K_m$, then 
$g(\node^{i+1}(\langle T,U\rangle))\ne g(\node^{i+1}(\langle T,U'\rangle))$ as they have different immediate successors. 
Similarly, if $r(T)\ne r(T')$ then $g(\node^{i+1}(\langle r(T),U\rangle))\ne g(\node^{i+1}(\langle r(T'),U'\rangle))$
as they have different immediate successors.   We can now conclude that $g$ preserves nonorder
by using the inductive hypothesis and the fact that $g$ preserves the immediate successor relation.
\end{proof}

\begin{lemma}\label{disjoint} The sets $\ran(g)$ and $k(\phi)$ are disjoint and $\ran(g)\cup k(\phi) = k(\psi)$.
\end{lemma}
\begin{proof}
It is immediate that $\ran(g)$ and $k(\phi)$ are disjoint.

We will first show that $\ran(g)\cup k(\phi)\subset k(\psi)$.  It is clear that $k(\phi)\subset k(\psi)$.
Since every node in $\ran(g)$ is at level $\geq i +1$, every node in $\ran(g)$ forces $\psi_1$.
Since $\psi_0$ is doubly negated and every successor of a node in $\ran(g)$
is in $\ran(g)$ or $k(\phi)$, by induction every element of $\ran(g)$ forces 
$\psi_0$ and $\ran(g)\subset k(\psi)$.

We will now show that $k(\psi)\subset \ran(g)\cup k(\phi)$. By construction, $k(\psi)\cap K^i_n
= k(\phi)\cup\{\gamma_A\}$ and $k(\psi)\cap \lev^n_{i+1} = \ran(g) \cap \lev^n_{i+1} = \ran(g|\lev^m_0) - \gamma_A$.
We will show that $k(\psi)\cap\lev^n_{j} \subset \ran(g|\lev^m_{j-(i+1)})$  for all $j \geq i + 1$ by induction 
on $j$.  We just observed that this holds for 
$j = i+1$. 

Suppose it holds for $j$. A node of $k(\psi)\cap\lev^n_{j+1}$ must be of the form 
$\node^{j+1}(\langle T,\emptyset\rangle)$ for $T\subset k(\psi)\cap K^j_n$. Since $T$ must contain an element
of $k(\psi)\cap\lev^n_j$ and every such node is below every element of 
$k(\psi)\cap\lev^n_{i-1}$, $T$ must be a subset of $k(\psi)\cap(K^j_n - 
K^{i-1}_n)$. Let $S = g^{-1}(T\cap(K^j_n - K^{i}_n))$ and $U = \{x_k \mid \alpha_k\in
T\cap \lev_i\}$.
Then $g^{-1}(\node^{j+1}_{K_n}(\langle T,\emptyset\rangle))$ is $\node^{j - i}_{K_m}(\langle S, \bigcap_{\mu\in S} w(\mu) - U\rangle)$

\end{proof}

It follows from Lemmas~\ref{injection} and~\ref{disjoint} that $g$ is an order-isomorphism from $K_m$ to $k(\phi) - k(\psi)$.
	
Define $f\colon F_m\to F_n$ by: 

1. $f(\bot) = \phi$ 

2. $f(x_i) = (\phi'_{\alpha_i}\vee\phi)\wedge\psi$

3. $f(\rho_0\wedge \rho_1) = f(\rho_0)\wedge f(\rho_1)$.

4. $f(\rho_0\vee \rho_1) = f(\rho_0)\vee f(\rho_1)$.

5. $f(\rho_0\rightarrow \rho_1) = (f(\rho_0)\rightarrow f(\rho_1))\wedge \psi$.

\begin{lemma}\label{fembedding} For any $\rho\in F_m$, $\phi\vdash f(\rho)\vdash \psi$.  If $\delta = g(\gamma)$ then
$\gamma \Vdash \rho$ iff $\delta \Vdash f(\rho)$.
\end{lemma}
\begin{proof}
The proof that $\phi\vdash f(\rho)\vdash \psi$ is an easy proof by induction on $\rho$.

We now prove the second part of the lemma by induction on $\rho$.

For $\rho = \bot$, the result is immediate. The observation that $\gamma\Vdash x_i$ iff
$\delta \not< \alpha_i$ furnishes the case where $\rho$ is $x_i$.  The inductive steps follow
from the existence of the order-isomorphism $g$ from $K_m$ to $k(\psi) - k(\phi)$ and
the fact that $\phi\vdash f(\rho)$ for all $\rho$.
\end{proof}

%
%
%
%

Note that it follows from Lemma~\ref{fembedding} that $f$ is injective and hence an embedding.

We now define a function from $F_n$ to $F_m$ that is
an inverse to $f$ when restricted to $[\phi,\psi]$. Define $h$ from $F_n$ to $F_{m}$ as follows:

1. $h(\bot) = h(x_i) = \bot$.

2. $h(\rho_0\wedge \rho_1) = h(\rho_0)\wedge h(\rho_1)$.

3. $h(\rho_0\vee \rho_1) = h(\rho_0)\vee h(\rho_1)$.


4. If there is some $\delta\in k(\phi)\cap K^{i-1}_n$ such that $\delta\not\Vdash \rho_0\rightarrow \rho_1$, then $h(\rho_0\rightarrow\rho_1) = \bot$. Otherwise, 
\[h(\rho_0\rightarrow\rho_1) =  (h(\rho_0)\rightarrow h(\rho_1))\wedge\bigwedge\{
x_i\mid \alpha_i\not\Vdash \rho\}\]

\begin{lemma}\label{inverse} Let $\delta = g(\gamma)$. For all $\rho\in F_n$, $\delta\Vdash \rho$ iff
$\gamma \Vdash h(\rho)$.
\end{lemma}
\begin{proof}
We will prove this by induction on the level of $\gamma$ and the structure of $\rho$.

If $\rho$ is $\bot$ or $x_i$, then $\delta\not\Vdash \rho$ and $\gamma\not\Vdash h(\rho)$.

The inductive step for $\rho = \rho_0\vee \rho_1$ and $\rho = \rho_0\wedge\rho_1$ is straightforward.

Let $\rho$ be $\rho_0 \rightarrow \rho_1$.  Suppose $\delta \Vdash \rho$. 
Then, since for every $\mu\in k(\phi)\cap K^{i-1}_n$, $\delta < \mu$, 
$h(\rho) = (h(\rho_0)\rightarrow h(\rho_1))\wedge\bigwedge\{
x_i\mid \alpha_i\not\Vdash \rho\}$. 
Since $\delta\Vdash \rho$, if $\alpha_i\not\Vdash \rho$, $\delta\not< \alpha_i$. It follows
that $\gamma\Vdash x_i$.  Thus $\gamma$ forces the right conjunct of $h(\rho)$.

Suppose $\delta \Vdash \rho_0$ and $\delta \Vdash \rho_1$. Then we are done by the inductive
hypothesis on the structure of $\rho$.  Otherwise, suppose $\delta\not\Vdash \rho_0$.
Then we are done by the inductive hypothesis on the structure of $\rho$ and the level
of $\gamma$.

Now suppose $\delta\not\Vdash\rho$.  Then there is some $\mu \geq \delta$
such that $\mu\Vdash \rho_0$ and $\mu\not\Vdash \rho_1$. If 
$\mu$ is in the range of $g$ then
we are done by induction. If $\mu\in K^{i-1}_n$, then $h(\rho) = \bot$ and we are done.
Otherwise $\mu\in K^i_n$ and is some $\alpha_j$.  Since $\delta < \alpha_j$, 
$\gamma \not\Vdash x_j$ and $\gamma\not\Vdash h(\rho)$.
\end{proof} 
It follows from Lemma~\ref{inverse} and Lemma~\ref{fembedding} that if $\phi\vdash \rho \vdash \psi$, then $f(h(\rho)) = \rho$.

\end{proof}
If $n \geq 2$, $m\geq 2$ by Fact~\ref{biglevel} we can find a level in $K_n$ satisfying the hypotheses of the
Proposition.  For example, we may pick a level in $K_n$ of cardinality greater than $2m$, call $2m$ of its
elements $\beta_1,\ldots, \beta_{2m}$, and let $A = \{\langle \{\beta_1,\beta_2\},\emptyset\rangle,
\ldots, \langle \{\beta_{2m - 1}, \beta_{2m}\},\emptyset\rangle\}$.

If $m = 1$, then we may let $\phi$ be $\bot$ and $\psi$ be $x_2\wedge \ldots\wedge x_n$.  The embedding $f$
from $H_1$ to $[\phi,\psi]\subset H_n$ sends $\rho$ to $\rho\wedge x_2\wedge\ldots \wedge x_n$. We may define a surjective
lattice homomorphism $h$ from $H_n$ to $H_1$ that is an inverse to $f$ as follows:

$h(x_1) = x_1$

$h(x_i) = \top$ for $1 < i \leq n$

$h(\phi\wedge \psi) = h(\phi)\wedge h(\psi)$

$h(\phi\vee \psi) = h(\phi)\vee h(\psi)$

$h(\phi\rightarrow \psi) = h(\phi)\rightarrow h(\psi)$
\end{proof}
\begin{corollary}\label{lattpreserve} There is a consequence-respecting, $(\wedge,\vee)$-preserving translation but not tautology-respecting
from $m$-variable intuitionistic logic to $n$-variable intuitionistic logic for
$n\geq 2$.
\end{corollary}
\begin{proof} Immediate.
\end{proof}
\begin{corollary}\label{tautpreserve} There is a consequence- and tautology-respecting translation
from $m$-variable intuitionistic logic to $n$-variable intuitionistic logic
for $n\geq 2$.
\end{corollary}
\begin{proof}
Let $f\colon F_m\to F_n$ be a consequence-respecting translation from
$m$-variable intuitionistic logic to $n$-variable intuitionistic logic.
Define $f'$ by $f'(\phi) = f(\top)\implies f(\phi)$. 

Then $f'$ is consequence- and tautology-respecting. To see that
it is consequence-respecting: If $\Gamma\vdash \phi$ then
$f(\Gamma)\vdash f(\phi)$, so $f(\top)\implies f(\Gamma),f(\top)\vdash f(\phi)$
and $f(\top)\implies f(\Gamma)\vdash f(\top)\implies f(\phi)$, where
$f(\top)\implies f(\Gamma)$ is an abbreviation of $\{f(\top)\implies \psi \mid
\psi\in f(\Gamma)\}$.

Conversely, if $f(\top)\implies f(\Gamma), f(\top)\vdash f(\phi)$,
then $f(\Gamma)\vdash f(\phi)$ since $f(\Gamma)\vdash f(\top)\implies f(\Gamma)$
and $f(\Gamma) \vdash f(\top)$ (this last fact is due to the fact that $f$ is 
consequence-preserving).
\end{proof}
\begin{corollary}\label{latttautpreserve} There is a consequence- and tautology-respecting,
$(\wedge,\vee)$-preserving translation from $m$-variable intuitionistic logic 
to $n$-variable intuitionistic logic for $n \geq 2$.
\end{corollary}
\begin{proof}
Let $f\colon F_m\to F_n$ be a consequence-respecting and $(\wedge,\vee)$-preserving
translation from $m$-variable intuitionistic logic to $n$-variable intuitionistic logic.
Define $f'$ by 
\[ f'(\phi) = \begin{cases} f(\phi) & \not\iv^m \phi \\
                            \top & \iv^m \phi
              \end{cases} \]
This is clearly still consequence-respecting and $\wedge$-preserving. The disjunction
property of intuitionistic logic implies that it is also $\vee$-preserving.
\end{proof}
Note that the translations given in Corollaries~\ref{lattpreserve} and~\ref{tautpreserve} can
be done in linear time, while the one given in Corollary~\ref{latttautpreserve} cannot, as it requires deciding whether
the given formula is a tautology. 

By~\cite{bellissima}, $H_n$ for $n\geq 2$ has an infinite descending chain, while $H_1$ does not, 
so there is no embedding of $H_n$ into $H_1$ for $n\geq 2$.


\section{A Lattice-Embedding from $H_\omega$ to $H_n$ for $n\geq 2$}
\begin{theorem}There is a lattice-embedding from $H_\omega$ into $H_2$
\end{theorem}
\begin{proof}
Pick $\alpha_1,\alpha_2,\alpha_3,\alpha_4,\alpha_5\in K_2$, all at the same level, say $i$. This may
be done by Fact~\ref{biglevel}.
Let $S = \{\rho \in K^i_2 \mid \forall i\in\{1,2,3,4\}\,\rho\not\geq \alpha_i\}$. Let
$\phi$ be
\[ (\neg\neg(\alpha_1\vee\alpha_2\vee\alpha_3\vee\alpha_4))\wedge \bigwedge_{\beta\in S} \phi'_\beta. \]
Let $T = \{\rho\in K^i_2 \mid \forall i\in \{1,2,3,4,5\}\,\rho\not\geq\alpha_i\}$. Let $\psi$
be
\[ (\neg\neg(\alpha_1\vee\alpha_2\vee\alpha_3\vee\alpha_4\vee\alpha_5))\wedge \bigwedge_{\beta\in T}\phi'_\beta.	 \]

Define a sequence $\{\beta^j_i\mid i\in \omega, j\in \{1,2,3,4\}\}$ as follows:
let $\beta^j_0 = \alpha_j$. For $i\geq 0$, let $\{\beta^j_{i+1}\mid j = 1,2,3,4\}$ be a collection
of four distinct nodes of the same level, with $\lev(\beta^1_{i+1})>\lev(\beta^1_i)$ and
such that they all force $\neg\neg(\beta^2_i\vee\beta^3_i\vee\beta^4_i)$. For example,
we may take $\beta^1_{i+1} = \langle \{\beta^2_i,\beta^3_i\},\emptyset\rangle$,
$\beta^2_{i+1} = \langle \{\beta^2_i,\beta^4_i\},\emptyset\rangle$,
$\beta^3_{i+1} = \langle \{\beta^3_i,\beta^4_i\},\emptyset\rangle$,
and $\beta^4_{i+1} = \langle \{\beta^2_i,\beta^3_i,\beta^4_i\},\emptyset\rangle$.

 As in~\cite{bellissima} (where a very similar construction is done),
the nodes of $\{\beta^1_i\mid i\in\omega\}$ are pairwise incomparable, and they all force $\phi$.

Define a Kripke model $K$ over the language $V_\omega = \{x_i \mid i \in \omega\}$ as follows: The set of nodes of $K$ is the set $k(\psi) - k(\phi)$
and a node $\alpha$ forces $x_i$ iff $\alpha\not\leq \beta^1_i$.

For all $\phi\in F_\omega$, let $k(\phi) = \{ \alpha\in K\mid \alpha\Vdash \phi\}$.

\begin{lemma} For all $\phi,\psi\in F_\omega$, $k(\phi)\subset k(\psi)$ iff $\phi\vdash \psi$.
\end{lemma}
\begin{proof}
Since $K$ is a Kripke model, if $\phi\vdash \psi$, $k(\phi)\subset k(\psi)$.

Suppose $\phi\not\vdash \psi$. Then there is a rooted finite Kripke model $K'$ over $V_\omega$ such that
$K' \Vdash \phi$ and $K'\not\Vdash \psi$. Since variables not occurring in $\phi$ or $\psi$
are irrelevant, we may assume that each node of $K'$ forces cofinitely many propositional variables.

Define a map $a\colon K'\to K$ inductively on $K'$ as follows: 
If $\gamma\in K'$ is a node such that $a(\gamma')$ has defined for all immediate successors
of $\gamma$, then let $a(\gamma)$ be a node whose set of successors
in $K_2$ is the upward-closure of the set 
$\{\beta^1_i \mid \gamma \not\Vdash x_i\}\cup \{a(\gamma') \mid \gamma'\geq \gamma\} \cup
\{\alpha_5\}$. 

For each $i$, $\gamma\Vdash x_i$ iff $a(\gamma)\Vdash x_i$.  Since $a$ is also order-preserving
and its range is upward-closed in $K$, we have that if $\gamma$ is the root of 
$K'$, $a(\gamma)\Vdash \phi$ and $a(\gamma)\not\Vdash \psi$.
\end{proof}

Now, as before, define
 $f\colon F_\omega\to F_2$ by: 

1. $f(\bot) = \phi$ 

2. $f(x_i) = (\phi'_{\beta^1_i}\vee\phi)\wedge\psi$

3. $f(\rho_0\wedge \rho_1) = f(\rho_0)\wedge f(\rho_1)$.

4. $f(\rho_0\vee \rho_1) = f(\rho_0)\vee f(\rho_1)$.

5. $f(\rho_0\rightarrow \rho_1) = (f(\rho_0)\rightarrow f(\rho_1))\wedge \psi$.

By precisely the same argument as before, this is an embedding.
\end{proof}

Note that, by~\cite{bellissima}, in any interval $[\phi,\psi]\subset H_n$,
there are atomic elements.  As there are no atomic elements in $H_\omega$,
$H_\omega$ cannot be embedded in $H_n$ as an interval.
%
%
%
%
%
%
%
%
%
%
%
\section{Impossibility of Lattice-Embedding $B_\omega$ into $H_\omega$}
Let $B_\omega$ be the countable atomless Boolean algebra. We will think of it as the Lindenbaum algebra
of classical propositional logic on countably infinitely many variables.
\begin{proposition}
There is no lattice embedding from $B_\omega$ into $H_n$ for any $n$ or into $H_\omega$.
\end{proposition}
\begin{proof}
By the previous theorem, it suffices to prove the proposition for $H_2$. Suppose there is a lattice embedding
of $B_\omega$ into $H_2$. Call it $f$.

Let $f(\top)$ have $n$ subformulas. Consider the $2^n$ formulas $\phi_1 = x_1\wedge\cdots\wedge x_n$,
$\phi_2 = x_1\wedge\cdots\wedge\neg x_n, \ldots, \phi_{2^n} = \neg x_1\wedge\cdots\wedge \neg x_n$.  
Since $f$ preserves $\wedge$ and $\vee$ we must have that $\{k(f(\phi_i))\mid 1 \leq i \leq 2^n\}$
is a partition of $k(f(\top)) - k(f(\bot))$ and that $k(f(\phi_i))\cap (k(f(\top)) - k(f(\bot)))$
is non-empty for each $i$.

\begin{lemma} Let $\alpha$ be a node in a Kripke model with exactly two immediate successors, $\alpha_1$
and $\alpha_2$.  Let $\phi$ be a formula.  Suppose that for each subformula $\phi'$ of $\phi$,
$\alpha_1\Vdash \phi'$ iff $\alpha_2\Vdash \phi'$ and that for all propositional variables $v$ appearing in $\phi$,
if $\alpha_1$ and $\alpha_2$ force $v$, then $\alpha\Vdash v$.  Then for each subformula $\phi'$ of
$\phi$, $\alpha \Vdash \phi'$ iff $\alpha_1\Vdash \phi'$.  In particular, $\alpha \Vdash \phi$
iff $\alpha_1\Vdash \phi$ iff $\alpha_2 \Vdash \phi$.
\end{lemma}
\begin{proof}
By induction on the structure of $\phi$. The conclusion is immediate if $\phi$ is atomic, and the 
$\wedge$ and $\vee$ cases are straightforward.

Suppose $\phi$ is $\phi_1\rightarrow \phi_2$. By induction, we can conclude that $\alpha\Vdash \phi'$
iff $\alpha_1\Vdash \phi'$ if $\phi'$ is a subformula of $\phi_1$ or $\phi_2$.  We just have to verify
that $\alpha\Vdash \phi$ iff $\alpha_1\Vdash \phi$.

Suppose $\alpha\Vdash \phi$.  Then, as $\alpha_1\geq \alpha$, $\alpha_1\Vdash \phi$.


Now suppose $\alpha\not\Vdash \phi$.  Thus, there must be some $\alpha' \geq \alpha$
such that $\alpha'\Vdash \phi_1$ and $\alpha'\not\Vdash \phi_2$. If $\alpha' = \alpha$ then
we are done by induction. Otherwise, we must have $\alpha' \geq \alpha_1$ or $\alpha' \geq \alpha_2$,
and thus $\alpha_1\not\Vdash \phi$.
\end{proof}
For each $i$, let $\beta_i\in k(f(\phi_i))\cap (k(f(\top)) - k(f(\bot)))$.  By the pigeonhole principle,
there must be some $i$ and $j$, $i\ne j$, such that $\beta_i\Vdash \phi'$ iff
$\beta_j\Vdash \phi'$ for all subformulas $\phi'$ of $f(\top)$. Let $\beta$ be 
$\langle \{\beta_i,\beta_j\},w(\beta_i)\cap w(\beta_j)\rangle$. We can easily verify that $\beta\in K_2$.
By the lemma, $\beta\in f(\top)$.  Thus, 
$\beta$ is in $k(f(\phi_m))\cap (k(f(\top)) - k(f(\bot)))$ for some $m$. Without loss
of generality, say $m \ne i$. Then $\beta_i \Vdash f(\phi_m)$ and $\beta_i \Vdash f(\phi_i)$
but $\beta_i \not\Vdash f(\bot)$, a contradiction.

\end{proof}
\section{Order-Embeddings}
\begin{proposition} Any countable partial ordering can be order-embedded into $H_2$ (and, therefore,
into $H_n$ for any $n\geq 2$).
\end{proposition}
\begin{proof}
We first make the following definition:
\begin{definition}[$\psi(\alpha_1,\ldots,\alpha_m)$, Permissive formulas] 
Let $\{\alpha_1,\ldots,\alpha_m\}$ be a set of nodes of $K_2$ all with the
same level. Let $S(\alpha_1,\ldots,\alpha_m) = 
 \{\delta \in K_2 \mid \lev(\delta)\leq \lev(\alpha_1)\text{ and }\forall i\, 
\delta\not\geq \alpha_i\}$.

We define $\psi(\alpha_1,\ldots,\alpha_m)$ to be
\[ \left(\neg\neg\bigvee_{i=1}^m\phi_{\alpha_i}\right)\wedge \bigwedge_{\delta\in S(\alpha_1,\ldots,\alpha_m)}\phi'_{\delta}\]
If $T = \{\alpha_1,\ldots,\alpha_m\}$ then $\psi(T)$ will denote
$\psi(\alpha_1,\ldots,\alpha_m)$. If some $\alpha_i$ is at a different
level than some $\alpha_j$, $\psi(\alpha_1,\ldots,\alpha_m)$ is not defined.

A formula of the form $\psi(\alpha_1,\ldots,\alpha_m)$ where $m\geq 3$ will be called \textit{permissive.}
The set $\{\alpha_1,\ldots,\alpha_m\}$ is called the set of generators of $\psi(\alpha_1,\ldots, \alpha_m)$
and $\lev(\alpha_1)$ is called the level of $\psi(\alpha_1,\ldots, \alpha_m)$.
\end{definition}

\begin{lemma} Given any permissive formula $\psi$, there exist permissive
formulas  $\psi_n$ for $n\in \{0,1\}$ such that for each $n\in\{0,1\}$, $k(\psi_n)\subset k(\psi)$ and 
$k(\psi_0)\cap k(\psi_{1})$ is finite.
\end{lemma}
\begin{proof}
Let $i$ be greater than the level of $\psi$ with $|\lev_i\cap k(\psi)| \geq 6$. We can find $i$
by Fact~\ref{bignotnot}.
Let $\lev_i\cap k(\psi) = \{\alpha_1,\alpha_2,\alpha_3,\beta_1,\beta_2,\beta_3,\ldots\}$,
let $\psi_0 = \psi(\alpha_1,\alpha_2,\alpha_3)$ and
$\psi_1 = \psi(\beta_1,\beta_2,\beta_3)$. 
\end{proof}
\begin{definition}[$\psi_\sigma$]
We define $\psi_\sigma$, for $\sigma\in \{0,1\}^{<\omega}$ as follows: 
Let $\psi_{\varepsilon} = \top$.
  Given $\phi_\sigma$, 
define $\phi_{\sigma n}$ for  $n\in\{0,1\}$ so that $\phi_{\sigma_n}$ is permissive, $k(\phi_{\sigma n}) \subset k(\phi_\sigma)$,
and $k(\phi_{\sigma 0})\cap k(\phi_{\sigma 1})$ is finite as in the above lemma.
\end{definition}
Note that $\psi_\sigma \vdash \psi_{\sigma'}$ iff $\sigma$ is an initial segment of $\sigma'$ as a binary string.
Note also that $\psi_\sigma \vdash \bigvee_i \psi_{\sigma_i}$ iff there is an $i$ such that
$\sigma = \sigma_i$.
\newcommand{\comp}[1]{\overline{#1}}
\begin{definition}[Complete Sets]
A set $S\subset H_2$ such that each element of $S$ is a disjunction  of the form $\bigvee_{i = 1}^n \psi_{\sigma_i}$ is called
complete if it satisfies the following property: Let $S_1$, $S_2$ be such that
$S_1\cup S_2 = S$, $S_1\cap S_2 = \emptyset$, $S_1$ is upward closed, and $S_2$ is downward closed.
Then 
there is some $\sigma(S_1,S_2)$ such that $|\sigma(S_1,S_2)| > \max\{|\sigma| \mid \psi_\sigma \text{ a
disjunct of a formula in } S\}$,
$\psi_{\sigma(S_1,S_2)}$ implies every element of $S_1$, and $\psi_{\sigma(S_1,S_2)}$ implies no element of
$S_2$.

Note that the condition that \[|\sigma(S_1,S_2)| > \max\{|\sigma| \mid \psi_\sigma \text{ a
disjunct of a formula in } S\}\] means
that no $\psi_\sigma\in S$ can imply $\psi_{\sigma(S_1,S_2)}$.

Note also that if $\sigma_1$ is such that $|\sigma_1| > \max\{|\sigma| \mid \psi_\sigma \text{ a
disjunct of a formula in } S\}$,
$\psi_{\sigma_1}$ implies every element of $S_1$, and $\psi_{\sigma_1}$ implies no element of
$S_2$ then so does $\sigma_1\sigma_2$ for any $\sigma_2$ (where the juxtaposition indicates
concatenation).  Without loss of generality, then, we may assume that each $\sigma(S_1,S_2)$
has the same length.
\end{definition}



The proposition will follow from the following lemma.
\begin{lemma} 
Suppose $P$ is a finite partial order, $S\subset H_2$ is a complete set, and $h$ is an isomorphism from $P$ to $(S,\leq)$. For any partial 
order $P'$ such that $|P'| = |P| + 1$, there is a $\phi$ of the form
$\bigvee_i \psi_{\sigma_i}$ such that $S\cup\{\phi\}$ is complete,
$P'\simeq (S\cup \{\phi\},\leq)$ via an isomorphism extending $h$.
\end{lemma}
\begin{proof}
If $P = \emptyset$, let $\phi$ be $\psi_0$. This is complete as we may let $\sigma(\{\phi\},\emptyset) = 00$
and $\sigma(\emptyset,\{\phi\}) = 10$.

%
%
Suppose $S = T_1\cup T_2 \cup T_3$ where $T_1$ is downward closed and $T_2$ is upward closed, and we would like to find $\phi$ so that $\phi$ is above all
the elements of $T_1$, below all the elements of $T_2$ and incomparable with the elements of $T_3$.


%
Let $\mathcal{S}$ be the collection of all partitions $(S_1,S_2)$ of $S$ such that $S_1$ is
upward closed and $S_2$ is downward closed.


Let $\phi$ be
\[ \bigvee_{\chi\in T_1} \chi \vee \bigvee \{\psi_{\sigma(S_1,S_2)0}\mid (S_1,S_2)\in \mathcal{S} \text{ and }
T_2\subset S_1\}\]


Clearly, $\phi$ is above every $\chi\in T_1$. We also have that $\phi$ is below every $\rho \in T_2$,
since every $\chi\in T_1$ must be below every $\rho \in T_2$, and by definition every $\psi_{\sigma(S_1,S_2)0}$
with $T_2\subset S_1$ is below every $\rho \in T_2$.


$\phi$ is not above any element in $T_2\cup T_3$: As noted above,
$\psi_\sigma \vdash \bigvee \psi_{\sigma_i}$ implies $\sigma = \sigma_i$ for some $i$.
But the disjuncts of $\phi$ are either elements of $T_1$ (which cannot be implied by elements of $T_2$ or $T_3$) or of the form $\psi_\sigma$ where
the length of $\sigma$ is greater than the length of any $\sigma'$ for $\psi_{\sigma'}$ some
disjunct of a formula in $T_2\cup T_3$.
%
%
%
%

$\phi$ is not below any element in $T_1\cup T_3$: Let $\mu \in T_1\cup T_3$. 
Let $S_2 = \{\mu' \mid \mu'\leq \mu\}$ and $S_1 = S- S_2$. Then $\psi_{\sigma'(S_1,S_2)0}$
is a disjunct of $\phi$ which does not imply $\mu$.

To see that $S\cup \{\phi\}$ is complete: Let $(S_1,S_2)\in \mathcal{S}$ with $\phi\in S_1$. 
Since $S_1$ is upward closed, we must have $T_2\cup S_1$.  Thus $\psi_{\sigma(S_1 - \{\phi\},S_2)0}$
is a disjunct of $\phi$. We may therefore take $\sigma(S_1, S_2)$ to be $\sigma(S_1 - \{\phi\}, S_2)0$
concatenated with enough zeroes to make its length greater than $\max\{\sigma \mid \psi_\sigma
\text{ a disjunct of a formula in } 
S\cup \{\phi\}\}$.


Let $(S_1, S_2)\in \mathcal{S}$ with $\phi\in S_2$. Thus $T_1\subset S_2$.
If there is any member of $T_2$ in $S_2$, then we may take
$\sigma(S_1,S_2)$ to be any sufficiently long extension
of $\sigma(S_1, S_2 - \{\phi\})$, since $\psi_{\sigma(S_1,S_2 - \{\phi\})}$
cannot imply $\phi$ since $\phi$ implies each element of $S_2$.

Thus we may assume that $T_2\subset S_1$. Therefore,
$\psi_{\sigma(S_1, S_2 - \{\phi\})0}$ is a disjunct of $\phi$. By construction,
there are no disjuncts of $\phi$ above it.  We may take $\sigma(S_1,S_2)$
to be any sufficiently long extension of $\sigma(S_1, S_2 - \{\phi\})1$.
%
\end{proof}

\end{proof}


\begin{thebibliography}{1}
\bibitem{coq}\textit{The Coq proof assistant}, \texttt{http://coq.inria.fr}.

\bibitem{nuprl}\textit{The PRL project}, \texttt{http://www.nuprl.org}.

\bibitem{beeson}Michael Beeson, \textit{Foundations of constructive mathematics: Metamathematical
studies}, Springer, Berlin/Heidelberg/New York, 1985

\bibitem{bellissima}Fabio Bellissima, \textit{Finitely generated free Heyting algebras}, Journal of Symbolic Logic \textbf{51} (1986), 152--165

\bibitem{butz}Carsten Butz, \textit{Finitely presented Heyting algebras}, 
\texttt{http://www.itu.dk/$\sim$butz/research/heyting.ps.gz}, 1998

\bibitem{darniere}Luck Darni\`ere and Markus Junker, \textit{On finitely generated Heyting
algebras}, \texttt{http://home.mathematik.uni-freiburg.de/junker/ preprints/heyting-221005.pdf}, 2005

\bibitem{ghilardi}Silvio Ghilardi and Marek Zawadowski, \textit{A sheaf representation and
duality for finitely presented Heyting algebras}, Journal of Symbolic Logic
\textbf{60} (1995), 911--939

\bibitem{girardpat}Jean-Yves Girard, Yves Lafont, and Paul Taylor, \textit{Proofs and types},
Cambridge University Press, Cambridge, 1989

\bibitem{nerodelambda}Anil Nerode, George Odifreddi, and Richard Platek, \textit{Constructive logics and lambda calculi}, in preparation

\bibitem{nishimura}Iwao Nishimura, \textit{On formulas of one variable in intuitionistic propositional calculus}, Journal of Symbolic Logic \textbf{25} (1960) 327--331

\bibitem{urq}Alasdair Urquhart, \textit{Free Heyting Algebras}, Algebra Universalis \textbf{3} (1973) 94--97
\end{thebibliography}
\end{document}